\documentclass[14pt,reqno]{article}
\usepackage{amsmath, amssymb}
\usepackage{latexsym, color}
\usepackage{epsfig}

\begin{document}

\newtheorem{thm}{Theorem}[section]
\newtheorem{cor}[thm]{Corollary}
\newtheorem{lem}[thm]{Lemma}
\newtheorem{prop}[thm]{Proposition}
\newtheorem{defn}[thm]{Definition}
\newtheorem{rem}[thm]{Remark}
\newtheorem{Ex}[thm]{EXAMPLE}
\def\nm{\noalign{\medskip}}

\bibliographystyle{plain}


\newcommand{\qed}{\hfill \ensuremath{\square}}
\newcommand{\ds}{\displaystyle}
\newcommand{\pf}{\medskip \noindent {\sl Proof}. ~ }
\newcommand{\p}{\partial}
\renewcommand{\a}{\alpha}
\newcommand{\z}{\zeta}
\newcommand{\pd}[2]{\frac {\p #1}{\p #2}}
\newcommand{\norm}[1]{\| #1 \|}
\newcommand{\dbar}{\overline \p}
\newcommand{\eqnref}[1]{(\ref {#1})}
\newcommand{\na}{\nabla}
\newcommand{\Om}{\Omega}
\newcommand{\ep}{\epsilon}
\newcommand{\tmu}{\widetilde \mu}
\newcommand{\vep}{\varepsilon}
\newcommand{\tlambda}{\widetilde \lambda}
\newcommand{\tnu}{\widetilde \nu}
\newcommand{\vp}{\varphi}
\newcommand{\RR}{\mathbb{R}}
\newcommand{\CC}{\mathbb{C}}
\newcommand{\NN}{\mathbb{N}}
\renewcommand{\div}{\mbox{div}~}
\newcommand{\bu}{{\bf u}}
\newcommand{\la}{\langle}
\newcommand{\ra}{\rangle}
\newcommand{\Scal}{\mathcal{S}}
\newcommand{\Lcal}{\mathcal{L}}
\newcommand{\Kcal}{\mathcal{K}}
\newcommand{\Dcal}{\mathcal{D}}
\newcommand{\tScal}{\widetilde{\mathcal{S}}}
\newcommand{\tKcal}{\widetilde{\mathcal{K}}}
\newcommand{\Pcal}{\mathcal{P}}
\newcommand{\Qcal}{\mathcal{Q}}
\newcommand{\id}{\mbox{Id}}
\newcommand{\be}{\begin{equation}}
\newcommand{\ee}{\end{equation}}

\title{Blow-up of Electric Fields between Closely Spaced Spherical Perfect Conductors}

\author{Mikyoung Lim\thanks{Department of Mathematics, Colorado State
University, Fort Collins, CO 80523, USA(lim@math.colostate.edu)}
\and KiHyun Yun\thanks{\footnotesize Department of Mathematics,
Michigan State University, East Lansing, MI 48824,
USA(kyun@math.msu.edu)}}

 \maketitle

\begin{abstract}
The electric field increases toward infinity in the narrow region
between closely adjacent perfect conductors as they approach each
other. Much attention has been devoted to the blow-up estimate,
especially in two dimensions, for the practical relevance to high
stress concentration in fiber-reinforced elastic composites. In this
paper, we establish optimal estimates for the electric field
associated with the distance between two spherical conductors in
$n-${ dimensional spaces for $n \geq 2$}. {The novelty of these
estimates is that they explicitly describe the dependency of the
blow-up rate on the geometric parameters: the radii of the
conductors.}
\end{abstract}

MSC-class: 15A15, 15A09, 15A23

\vspace{2pc}

\section{Introduction}
We consider the blow-up of the electric fields in the narrow region
between a pair of perfect conductors which is closely adjacent in
$n$ dimensions ($n\geq 2 $). Conductors provide higher intensity of
electric flux around them. The intensity increases as a pair of
conductors approaches each other, and the electric field even
reaches toward infinity (refer to \cite{AKL, AKLLL,Y,
Y2,ADKL,BLY})).

\par In this paper, we present the optimal blow-up estimate for the
electric field with respect to the distance between a pair of
conductors under the assumption that the conductors are of spherical
shape in $n$ dimensions ($n\geq 2$). {The novelty of these estimate
is to describe explicitly the dependency of the blow-up rate on the
radii of the conductors: this paper is the first result to establish
the role of the geometrical factor of conductors in the blow-up of
the electric field in three or higher dimensions.}

\par {Besides the consideration of the gradient estimates in the frame of the electrostatic theory, much attention has been drawn to it of the relevance to
the stress--strain behavior of composite materials, especially in two dimensions.}
 According to Budiansky and Carrier \cite{BC}, unexpectedly low strengths in stiff fiber-reinforced
 composites have been reported, due to the high stress concentration occurring
 in the narrow region between fibers (also refer to \cite{K}). In the anti-plane shear model,
 the stress tensor represents the electric field in the two dimensional conductivity model, where the out-of-plane elastic displacement satisfies a conductivity
equation  \cite{BV}. Thus, {the gradient estimates for electric
field have a valuable meaning in relation to in the failure analysis
of composite material.} To give a brief description of related
works, for the case that the inclusions and the outside of
inclusions have the comparable conductivities (or shear moduli), {
it was verified} that the electric field remains bounded
independently of the distance between the inclusions. Li and
Vogelius \cite{LV} have shown that the electric field does not blow
up even when the inclusions approaches each other. Moreover, Li and
Nirenberg \cite{LN} have extended this result to elliptic systems.
These results point out that the extremely high conductivity (or the
stiffness of fibers) is indispensable to the blow-up phenomena.

\par In this respect, much attention has been focused on the model of
a pair of perfect conductors which are $\epsilon$ apart. Ammari,
Kang, H. Lee, J. Lee and Lim \cite{AKL,AKLLL} have established the
optimal blow-up rate $\epsilon^{-1/2}$ as the distance $\epsilon$
goes to zero, when conductors are of circular shape in two
dimensions. Yun \cite{Y,Y2} has extended the above mentioned result
to a sufficiently general class of the conductors' shapes in two
dimensions. In three or higher dimensional case, Bao, Li and Yin
\cite{BLY} recently obtained the optimal blow-up rate for perfect
conductors of general shape: the optimal blow-up rate is
$\left(\epsilon |\log \epsilon|\right)^{-1}$ for three dimensions,
and is $\epsilon^{-1}$ for higher $n$ dimensions ($n\geq 4$).
However, their estimates are only given by the distance between two
conductors and geometric information of conductors are not
incorporated into the estimates.

\par{What is new in this paper is that for the case of spherical perfect conductors in three and higher dimensions, the gradient estimates
are established in terms of the radii as well as the distance between inclusions.
What is more is that the approach introduced in this paper to derive the estimates is distinct from
the methods of Bao et al. \cite{BLY} and Ammari et al. \cite{AKL, AKLLL}.
In the two dimensional case, our approach provides the same estimates as of Ammari et el., Proposition \ref{p32}, in a much simpler way for the case of perfect inclusions.
}

\section{Mathematical formalism and main results}
\par From now on, $\mathbb{R}^n$  denotes  $n$ dimensions, and $B_r (x_1, x_2, \cdots, x_n)$
 { is} the sphere with
radius $r$ and center $(x_1, x_2, \cdots, x_n)$ in $\mathbb{R}^n$.
Given any entire harmonic function $H$ in $\mathbb{R}^n~(n\geq 2)$,
we define the electric potential $u$ as the unique solution {to the
following conductivity problem}:
\begin{equation} \label{eq:001}
\quad \left\{
\begin{array}{ll}
\ds\Delta u  = 0,\quad&\mbox{in }{\mathbb{R}^{n} \backslash \overline{(D_1 \cup D_2)}},\\
\ds u(\textbf{x})- H(\textbf{x}) = O(|\textbf{x}|^{1-n})\quad&\mbox{as } |\textbf{x}| \rightarrow \infty,\\
\ds u |_{\partial D_i} = C_i~ \mbox{(constant)},\\
 \int_{\partial D_i} {\partial_{  \nu} u }~dS = 0,\quad
&\mbox{for } i = 1, 2,
\end{array}
\right.
\end{equation} where $\textbf{x} =
(x_1,x_2,\cdots, x_n)$. This solution $u$ can be interpreted
physically as the electric potential outside conductors $D_1$ and
$D_2$ under the action of applied electric field $\nabla H$.

\par In this paper, we start by considering the case {that} $\nabla H$ is a
uniform field, i.e, $H= \mathbf{a} \cdot \textbf{x}$ for some
constant $\mathbf{a}$ in $\mathbb{R}^n${, in  Theorem \ref{thm2} and
\ref{thm3}. Based on these, the optimal upper bound of the gradient
for any entire harmonic function $H$ is established in Theorem
\ref{thm4}.}

\begin{thm}[Three dimensions]\label{thm2} We assume that $D_1$ and $D_2$ are the pair of spheres with radii $r_1$ and
 $r_2$ that are $2 \epsilon$ apart in $\mathbb{R}^3$. Thus, we set
  $$D_1 = B_{r_1}(r_1 + \epsilon, 0,0)~\mbox{and}~D_2 = B_{r_2}(-(r_2+  \epsilon), 0,0).$$
  Let $u$ be the solution to (\ref{eq:001}) for $H(x_1,x_2,x_3) = \sum_{i=1}^{3} a_i x_i$.
Then, there exists a { positive} constant $C_*$ independent of
$\epsilon$, $r_1$, $r_2$ and $(a_1,a_2,a_3)$ such that
$$ \frac{1}{C_*}|a_1|\left(\frac {r_1 r_2}{r_1 + r_2} \right)\frac{1}{| \log \epsilon |}\leq\Big| u|_{\partial D_1 } - u|_{\partial D_2}
\Big| \leq {C_*}|a_1|\left(\frac {r_1 r_2}{r_1 + r_2} \right) {|
\log \epsilon|}$$ for sufficiently small $\epsilon > 0$ . {
\begin{enumerate}
\item[(a)]
\par In the case that $a_1$ is nonzero, for any sufficiently small
$\epsilon$, there is a point $\mathbf{x}_0$ between $D_1$ and $D_2$
such that
$$\frac 1 {2 C_{*}} |a_1|\left(\frac {r_1 r_2}{r_1 + r_2} \right) {\frac 1
{|\epsilon \log \epsilon |}} \leq | {\nabla u (\mathbf{x}_0) }|.$$
The lower bound above is optimal in the sense that there is a positive constant $C^*$ independent of $\epsilon$, $r_1$, $r_2$ and
$(a_1, a_2, a_3)$, satisfying that
 \be \norm{\nabla u }_{L^{\infty}
(\mathbb{R}^{3}\setminus (D_1 \cup D_2))} \leq C^* |a_1|\left(\frac
{r_1 r_2}{r_1 + r_2} \right) {\frac 1 {|\epsilon \log \epsilon |}},
\label{main}\ee for sufficiently small $\epsilon > 0$.
\item[(b)]
In the case that $a_1$ is zero, the gradient of $u$ does
not blow up even when the distance $\epsilon$ goes to zero, i.e.,
there is a positive constant $C_0 ^*$ independent of
$\epsilon$ and $(a_1, a_2, a_3)$, satisfying that \be  \norm{\nabla u
}_{L^{\infty} (\mathbb{R}^{3}\setminus (D_1 \cup D_2))} \leq C_0 ^*
\left(|a_2|+|a_3|\right)\label{eq:110},\ee for sufficiently small
$\epsilon
> 0$ .
\end{enumerate}}
\end{thm}
\begin{rem}\label{rem}
The constant $C_0 ^*$ at \eqref{eq:110} depends on $r_1$ and $r_2$:
in details, there is a constant $C$ so that \be\norm {\nabla u
}_{L^{\infty} (\mathbb{R}^3 \setminus (D_1 \cup D_2))} \leq C \max
\left\{\frac {r_1}{r_2}, \frac {r_2}{r_1}  \right\}{
\left(|a_2|+|a_3|\right)}, \label{eq:111}\ee when $a_1=0$. The
derivation of the inequality above is included in the proof of
Theorem \ref{thm2}. The term '$C \max \left\{\frac {r_1}{r_2}, \frac
{r_2}{r_1}  \right\}$' becomes arbitrarily large for small $r_1$ or
$r_2$, and { then it is not guaranteed that} the bound
\eqref{eq:111} above is optimal. However, our attention is focused
on the contribution of $r_1$ and $r_2$ to the blow-up rate
associated to the distance $2\epsilon$. As mentioned in Theorem
\ref{thm2}, the gradient of $u$ blows up as $\epsilon \rightarrow 0$
if and only if $a_1$ is nonzero. On this occasion, the $ 1
/{|\epsilon \log \epsilon|}$ term at \eqref{main} dominates the
upper bound of $\norm{\nabla u }_{L^{\infty}
(\mathbb{R}^{3}\setminus (D_1 \cup D_2))}$, and the bound
\eqref{main} describes the contribution of $r_1$ and $r_2$ well.
Therefore, Theorem \ref{thm2} would suffices for our purpose.

\end{rem}

\begin{thm}[Higher dimensions]\label{thm3} We assume that $D_1$ and $D_2$ are the pair of spheres with radii $r_1$ and
 $r_2$ that are $2 \epsilon$ apart in $\mathbb{R}^n~(n\geq 4)$. Thus, we set
  $$D_1 = B_{r_1}(r_1 + \epsilon, 0,\cdots,0)~\mbox{and}~D_2 = B_{r_2}(-(r_2+  \epsilon), 0,\cdots,0).$$
  Let $u$ be the solution to (\ref{eq:001}) for $H(x_1,x_2,\cdots,x_n) =\sum_{i=1}^{n} a_i x_i$.
Then, there exists a constant $C_{**}$ independent of $\epsilon$,
$r_1$, $r_2$ and $(a_1,a_2,\cdots,a_n)$ such that
$$ \frac{1}{C_{**}}|a_1|\left(\frac {r_1 r_2}{r_1 + r_2} \right)\leq\Big| u|_{\partial D_1 } - u|_{\partial D_2}
\Big| \leq {C_{**}}|a_1|\left(\frac {r_1 r_2}{r_1 + r_2} \right) {
}$$ for sufficiently small $\epsilon > 0$ .

\begin{enumerate}
{
\item[(a)]
\par In the case that $a_1$ is nonzero, for any sufficiently small
$\epsilon$, there is a point $\mathbf{x}_0$ between $D_1$ and $D_2$
such that
$$\frac 1 {2 C_{**}} |a_1|\left(\frac {r_1 r_2}{r_1 + r_2} \right) {\frac 1
{\epsilon }} \leq | {\nabla u (\mathbf{x}_0) }|.$$  The lower bound above is optimal in the sense that there is a positive
constant $C^{**}$ independent of $\epsilon$, $r_1$, $r_2$ and
$(a_1,a_2,\cdots,a_n)$, satisfying that
 $$ \norm{\nabla u }_{L^{\infty}
(\mathbb{R}^{3}\setminus (D_1 \cup D_2))} \leq C^{**}
|a_1|\left(\frac {r_1 r_2}{r_1 + r_2} \right) {\frac 1 {\epsilon
}},$$ for sufficiently small $\epsilon > 0$.

\item[(b)]
In the case that $a_1$ is zero, the gradient of $u$ does
not blow up even when the distance $\epsilon$ goes to zero, i.e., there is a positive constant $C_0 ^{**}$ independent of
$\epsilon$, $r_1$, $r_2$ and $(a_1,a_2,\cdots,a_n)$, satisfying
$$  \norm{\nabla u }_{L^{\infty}
(\mathbb{R}^{3}\setminus (D_1 \cup D_2))} \leq C_0 ^{**} \sum_{i=2}
^n|a_i|,$$ for sufficiently small $\epsilon > 0$ .
}
\end{enumerate}
\end{thm}
Similarly to Remark \ref{rem}, the constant $C_0 ^{**}$ above
depends on $r_1$ and $r_2$,i,e, there is a constant $C$ so that
\be\norm {\nabla u }_{L^{\infty} (\mathbb{R}^n \setminus (D_1 \cup
D_2))} \leq C \max \left\{\frac {r_1}{r_2}, \frac {r_2}{r_1}
\right\} \sum_{i=2} ^{n}| a_i|, \notag\ee when $a_1=0$. The
derivation of the inequality above is also presented in the proof of
Theorem \ref{thm3}.

\par  In this work, the blow-up estimates in terms of radii is
presented only for
 three or higher dimensional case. Speaking of two dimensions, Ammari \textit{et al}
\cite{AKL} already provided the optimal bound  (\ref{eq:109}) and
 (\ref{eq:108}) in terms of radii of circular inclusions as
\be\max |\nabla (u - H) | \leq C \sqrt {\frac {r_1 r_2}{ r_1 + r_2
}} \frac {1} {\sqrt {\epsilon }} \label{eq:109}\ee where $r_1$ and
$r_2$ are the radii of circular inclusions. This is also derived in
Proposition \ref{p32} of this paper. As has been mentioned before,
the method in this paper is much simpler method.

\begin{thm}[General entire harmonic function $H$] \label{thm4} Let $D_1$ and $D_2$ be a pair of balls as assumed in the previous
theorems in $n$ dimensions ($n\geq 3$). We choose a large bounded
domain $\Omega$ containing $D_1$ and $D_2$ for any small
$\epsilon>0$. For the sake of convenience, we select the ball
$B_{4(r_1 + r_2)}(0,\cdots,0)$ as $\Omega$. For any given entire
harmonic function $H$, let $u$ be the solution to (\ref{eq:001}) for
$H$.
\par In three dimensions $(n=3)$, there is a constant $C^{*}$, independent of $\epsilon$, $r_1$, $r_2$ and
$(a_1,a_2,a_3)$, satisfying
 $$ \norm{\nabla (u -H  ) }_{L^{\infty}
(\mathbb{R}^{3}\setminus (D_1 \cup D_2))} \leq C^{*} \norm {\nabla
H}_{L^{\infty} (\Omega)}\left(\frac {r_1 r_2}{r_1 + r_2} \right)
{\frac 1 {|\epsilon \log \epsilon |}}$$ and
$$\norm{\nabla u }_{L^{\infty}
(\Omega\setminus (D_1 \cup D_2))} \leq C^{*}\norm {\nabla
H}_{L^{\infty} (\Omega)}\left(\frac {r_1 r_2}{r_1 + r_2} \right)
{\frac 1 {|\epsilon \log \epsilon |}}$$ for sufficiently small
$\epsilon > 0$.
\par In higher dimensions $(n\geq 4)$, there is a constant $C^{**}$, independent of $\epsilon$, $r_1$, $r_2$ and
$(a_1,a_2,\cdots,a_n)$, satisfying
 $$ \norm{\nabla (u -H  ) }_{L^{\infty}
(\mathbb{R}^{3}\setminus (D_1 \cup D_2))} \leq C^{**} \norm {\nabla
H}_{L^{\infty} (\Omega)}\left(\frac {r_1 r_2}{r_1 + r_2} \right)
{\frac 1 {\epsilon}}$$ and
$$\norm{\nabla u }_{L^{\infty}
(\Omega\setminus (D_1 \cup D_2))} \leq C^{**} \norm {\nabla
H}_{L^{\infty} (\Omega)}\left(\frac {r_1 r_2}{r_1 + r_2} \right)
{\frac 1 {\epsilon }}$$ for sufficiently small $\epsilon > 0$.

\end{thm}
\section{Representation of the potential difference}
We introduce a harmonic function $h$ as follows:
\begin{equation} \label{eq:101}
\quad \left\{
\begin{array}{ll}
\ds\Delta h  = 0,\quad&\mbox{in }{\mathbb{R}^n \backslash \overline{(D_1 \cup D_2)}}, \\
\ds h= O(|\textbf{x}|^{1-n}),\quad&\mbox{as } |\textbf{x}| \rightarrow \infty,\\
\ds h |_{\partial D_i} = k_i\mbox{ (constant)},\\
\int_{\partial D_i} {\partial_{  \nu} h }~dS =
(-1)^{i+1},\quad&\mbox{for }  i = 1, 2.
\end{array}
\right.
\end{equation}
It is essential in this work to construct the function $h$ because
of the following lemma.
\begin{lem}\label{lemma:h}(\cite{Y})
For a solution $u$ to \eqnref{eq:001}, we have that
\begin{equation}\label{eq:102}
u\big|_{\partial D_1} - u\big|_{\partial D_2}=\int_{\partial D_1}
({\partial_{  \nu} h }) H\ dS+ \int_{\partial D_2}( {\partial_{ \nu}
h }) H\ dS.\end{equation}
\end{lem}
\pf With the boundary condition of $u$ on $\p D_1$ and $\p D_2$, Green's
identity for $H$ inside $D_1$ and $D_2$ yields
$$\int_{\partial D_1 } \partial_{\nu} (u-H)\  dS =\int_{\partial D_2 } \partial_{\nu} (u-H)\ dS=0,$$
and thus
$$I:=\int_{\partial D_1 } h\partial_{\nu} (u-H)\ dS +\int_{\partial D_2 }h \partial_{\nu} (u-H)\ dS=0.$$

Applying again Green's identity outside $\overline{D_1\cup D_2}$,
we have
\begin{align*}
0&=I=\big(u\big|_{\partial D_1} - u\big|_{\partial
D_2}\big)-\int_{\partial D_1} ({\partial_{  \nu} h }) H\ dS-
\int_{\partial D_2}( {\partial_{ \nu} h }) H\ dS.
\end{align*}\qed

We remark that the above representation \eqnref{eq:102} is
observed by Yun \cite{Y} for the purpose of estimating the
stresses between two arbitrary shaped inclusions in $\mathbb{R}^2$. By
constructing a harmonic function $h$ and calculating the right
hand side of \eqnref{eq:102}, Yun estimated the potential
difference between two adjacent conductors.

The idea to establish $h$ is from the basic theory in
electrodynamics, and we use several times the following property
of Apollonius circles.
\subsection{Apollonius Circle in $\RR^n$}
For a ball $B_r(\mathbf{c})$ in $\RR^n$ and a point
$\mathbf{p}$, $|\mathbf{p}-\mathbf{c}|>r$, we have
\begin{equation}\label{apoll}
\frac{r}{|\mathbf{p}-\mathbf{c}|}\frac{1}{|\mathbf{x}-R(\mathbf{p})|}=\frac{1}{|\mathbf{x}-\mathbf{p}|},
\quad \mbox{for all }\mathbf{x}\in\p B_r(\mathbf{c}),\end{equation}
where
$R$ is the reflection with respect to $B_r(\mathbf{c})$, i.e.,
\begin{equation*}
R(\mathbf{p})=\frac{r^2(\mathbf{p}-\textbf{c})}{|\mathbf{p}-\textbf{c}|^2}+
\mathbf{c}.
\end{equation*}

%

A simple application of Apollonius circle is {estimating} the
potential difference of the solution to \eqnref{eq:001} for two
circles with different radii.

\subsection{Estimates in $\mathbb{R}^2$}
We let $$D_1 = B_{r_1}(\mathbf{c}_1)~\mbox{and}~D_2 =
B_{r_2}(\mathbf{c}_2),$$
where $\mathbf{c}_1=(r_1 + \epsilon, 0)$ and $\mathbf{c}_2=(-r_2- \epsilon,0)$, and $R_i$ be the reflection with respect to $D_i$, in other words,
$$R_i(\mathbf{x})=\frac{r_i^2(\mathbf{x}-\mathbf{c}_i)}{|\mathbf{x}-\mathbf{c}_i|^2}+ \mathbf{c}_i,\ i=1,2. $$

Let $\textbf{p}_1\in D_1$ be the fixed point of $R_1 \circ R_2$, then
$R_2(\textbf{p}_1)(=:\textbf{p}_2)$ is the fixed point of  $R_2 \circ R_1$ and $R_1(\textbf{p}_2)=\textbf{p}_1$ .
From \eqnref{apoll}, \begin{equation*}
\frac{|\mathbf{x}-\mathbf{p}_1|}{|\mathbf{x}-\mathbf{p}_2|}=\left\{
\begin{array}{ll}
\ds\frac{r_1}{|\mathbf{p}_2-\mathbf{c}_1|},&\mbox{ for }\mathbf{x}\in\p D_1,\\[6mm]
\ds\frac{|\mathbf{p}_1-\mathbf{c}_2|}{r_2},&\mbox{ for }\mathbf{x}\in\p D_2.
\end{array}\right.
\end{equation*}
Hence, the solution to \eqnref{eq:101} is
$$h:=\frac 1 {2\pi} \left( \log |\textbf{x}-\textbf{p}_1|- \log |\textbf{x}-\textbf{p}_2|\right)
=\frac{1}{2\pi}\log\left(\frac{| \textbf{x}-\textbf{p}_1|}{|\textbf{x}-\textbf{p}_2|}\right),$$
and, from \eqnref{eq:102}, we have the following proposition.
\begin{prop}\label{p32}
Let $H(x_1,x_2)$ be an entire harmonic function. The solution $u$ to
(\ref{eq:001})
 satisfies
 \be
 \begin{array}{ll}u|_{\partial D_1} - u|_{\partial D_2} &=
H(\textbf{p}_1)-H(-\textbf{p}_2)\\ &= 4 \partial_{x_1} H (0,0)
\sqrt{\frac{r_1r_2}{r_1+r_2}}\sqrt\epsilon+O(\epsilon).
 \label{eq:107}\end{array}\ee
 Referring to the mean value theorem, there exists a point
$\textbf{x}_2$ between $\partial D_1$ and $\partial D_2$ such that
\be |\nabla u (\textbf{x}_2)|\geq |\partial_{x_1} H(0,0)|
\sqrt{\frac{r_1r_2}{r_1+r_2}} \frac 1 {\sqrt{\epsilon}}.
\label{eq:108}\ee for any sufficiently small $\epsilon >0$.
Moreover, there is a constat $C$ independent of $\epsilon$, $r_1$
and $r_2$ such that
$$ \norm{\nabla u }_{L^{\infty}(\Omega \setminus (D_1 \cup D_2)) }\leq C \norm {\nabla H}_{L^{\infty}(\Omega)}
\sqrt{\frac{r_1r_2}{r_1+r_2}} \frac 1 {\sqrt{\epsilon}} $$ where
$\Omega=B_{4(r_1+r_2)}(0,0).$
\end{prop}
\pf
The fixed points $\textbf{p}_i$ satisfies
$$\textbf{p}_1 =
\Bigr(2\sqrt{\frac{r_1r_2}{r_1+r_2}}\sqrt\epsilon +O(\epsilon)
,0\Bigr)\mbox{ and }\textbf{p}_2 =
\Bigr(-2\sqrt{\frac{r_1r_2}{r_1+r_2}}\sqrt\epsilon +O(\epsilon)
,0\Bigr).$$ Therefore, we obtain (\ref{eq:107}) and (\ref{eq:108}).
By virtue of the argument presented by Bao et al. in \cite{BLY}, the
upper bound of the gradient is derived from (\ref{eq:107}). In this
paper, the same process as this proposition to derive the upper
bound of the gradient is also presented in the proof of Theorem
\ref{thm2}. Therefore, please refer to the derivation of the upper
bound of the gradient in the proof of Theorem  \ref{thm4}. \qed

We remark that the same gradient estimate has been obtained by
Ammari et al. in \cite{AKL}. They represented $u$ by single layer
potentials of the Laplacian with potential functions defined using
Kelvin transform $R_i$, $i=1,2$, and obtained the blow-up rate by
investigating the potential functions. The novelty of their work is
that their estimates is not only for a extreme conductivity but also
for a finite positive constant. However, as for the extreme case,
our result provides a much simpler method for obtaining the blow-up
rate.

\section {Derivation for Theorem \ref{thm2}}\label{sec:b}
Differently from the two dimensional space where the point charge
potential, the logarithm, separate the multiplication with ratio
$\rho$ into a sum, we cannot constructed $h$ just with two point
charge potential functions in higher dimensional space. Therefore, we
introduce a sequential process to build $h$.

\subsection{Construct $h$ in $\RR^n$, $n\geq3$}\label{subsec:h}
Let $$D_1 = B_{r_1}(\mathbf{c}_{1})~\mbox{and}~D_2 =
B_{r_2}(\mathbf{c}_{2}),$$
 where $$\mathbf{c}_1=(r_1+\epsilon,0,\dots,0)\mbox{ and }\mathbf{c}_2=(-r_2-\epsilon,0,\dots,0).$$

 We start from a  harmonic function $h_{1,0}$, defined outside of $\bar{D}_1$,  which is $$h_{1,0}(\textbf{x})={\frac {1}{|\textbf{x}-\mathbf{c}_{1}|^{n-2}}}.$$
 Note that $h_{1,0}$ is constant on $\p D_1$. However, it  is not constant on
$\p D_2$, and  we neutralize it by adding auxiliary point charge potential
$h_{1,1}$ to make ($h_{1,0}+h_{1,1}$) constantly zero on $\p
D_2$. From \eqnref{apoll}, $h_{1,1}$ is defined as
$$h_{1,1}(\textbf{x})=\left(\frac{r_2}{|\mathbf{c}_2-\mathbf{c}_{1}|}\right)^{n-2}\frac{-1}{|\textbf{x}-R_2(\mathbf{c}_1)|^{n-2}},$$
where  $R_i$, $i=1,2$, be the reflection with respect to $D_i$. On
the next step, we add {$h_{1,2}$ to ($h_{1,0}+h_{1,1}$) and make
($h_{1,0}+h_{1,1}+h_{1,2}$) be constant on} $\p D_1$, i.e.,
$$h_{1,2}=\left(\frac{r_2}{|\mathbf{c}_2-\mathbf{c}_{1}|}\right)^{n-2}\left(\frac{r_1}{|\mathbf{c}_1-R_2(\mathbf{c}_1)|}\right)^{n-2}
\frac{1}{|\textbf{x}-R_1(R_2(\mathbf{c}_1))|^{n-2}}$$
Consequently, we construct $h_{1,m}$ , $m\in\NN$ as
\begin{equation}\label{eqn:h1m}
h_{1,m}=\bigr(q_{1,m}\bigr)^{n-2} \frac{(-1)^m}{|\mathbf{x}-\mathbf{c}_{1,m}|^{n-2}},
\end{equation}
where
\be\label{c1m}
\mathbf{c}_{1,m}=
\begin{cases}
(R_1R_2)^k(\mathbf{c}_1),&\quad\mbox{if }m=2k, \ k\geq0,\\
R_2(R_1R_2)^{k}(\mathbf{c}_1),&\quad \mbox{if }m=2k+1,\ k\geq0,
\end{cases}
\ee
%
\begin{equation}\label{q1m}
q_{1,m}=\prod_{j=0}^{m}\rho_{1,j},\quad \mbox{for }m\in\NN,
\end{equation}
and
\be\label{rho1}
\ds
\rho_{1,j}=
\begin{cases}
\ds 1,&\quad\mbox{if }j=0,\\
\ds\frac{r_1}{|\mathbf{c}_1-\mathbf{c}_{1,2k-1}|},&\quad\mbox{if }j=2k,\ k\geq1\\
\ds\frac{r_2}{|\mathbf{c}_2-\mathbf{c}_{1, 2k}|},&\quad\mbox{if }j=2k+1\ k\geq0.
\end{cases}
\ee

Similarly, we define $h_{2,m}$ , $m\in\NN$ as
\begin{equation}\label{eqn:h2m}
h_{2,m}=\bigr(q_{2,m}\bigr)^{n-2} \frac{(-1)^m}{|\mathbf{x}-\mathbf{c}_{2,m}|^{n-2}},
\end{equation}
where
\be\label{c2m}
\mathbf{c}_{2,m}=
\begin{cases}
(R_2R_1)^k(\mathbf{c}_2),&\quad\mbox{if }m=2k, \ k\geq0,\\
R_1(R_2R_1)^{k}(\mathbf{c}_2),&\quad \mbox{if }m=2k+1,\ k\geq0,
\end{cases}
\ee
%
\begin{equation}\label{q2m}
q_{2,m}=\prod_{j=0}^{m}{\rho_{2,j}},\quad \mbox{for }m\in\NN,
\end{equation}
and
\be\label{rho2}
\ds
\rho_{2,j}=
\begin{cases}
\ds 1,&\quad\mbox{if }j=0,\\
\ds\frac{r_2}{|\mathbf{c}_2-\mathbf{c}_{2,2k-1}|},&\quad\mbox{if }j=2k,\ k\geq1\\
\ds\frac{r_1}{|\mathbf{c}_1-\mathbf{c}_{2, 2k}|},&\quad\mbox{if }j=2k+1\ k\geq0.
\end{cases}
\ee

\smallskip

Since $(R_1R_2)^k(\mathbf{c}_{1})\in D_1$ and $R_2(R_1R_2)^k(\mathbf{c}_{1})\in D_2$, we have, for $j\geq1$,
$$\rho_{1,j}\leq  \max_{i=1,2}\frac{r_i}{r_i+2\epsilon}=\frac{1}{1+\frac{2\epsilon}{r_{\max}}},\mbox{ where }r_{\max}=\max(r_1,r_2).$$
By the same way, $\rho_{2,j}\leq\frac{1}{1+\frac{2\epsilon}{r_{\max}}}$. Hence,
$$\sum_{m=0}^{\infty}\bigr(q_{s,m}\bigr)^{n-2}<\infty,\quad s=1,2$$
and the two series $\sum_{m=0}^\infty h_{i,m},\ i=1,2$, are well
defined. To get $h$ satisfying the decaying condition, we sum the
series $h_{1,m}$ and $h_{2,m}$ with different weights as the
following lemma.
\begin{lem}\label{series:h}
The solution to \eqnref{eq:101} is given by
\begin{equation}
\ds h(\mathbf{x})=\frac{1}{(2-n)\omega_n}\frac{1}{ M}
\left[Q_2\sum_{m=0}^\infty \frac{(-1)^m\bigr(q_{1,m}\bigr)^{n-2} }{|\mathbf{x}-\mathbf{c}_{1,m}|^{n-2}}
-Q_1\sum_{m=0}^\infty \frac{(-1)^m\bigr(q_{2,m}\bigr)^{n-2}}{|\mathbf{x}-\mathbf{c}_{2,m}|^{n-2}}\right],
\end{equation}
where $\omega_n$ is the area of the unit sphere, and
\begin{equation*}\label{def:Q}
Q_s=\sum_{m=0}^{\infty}\Bigr[(-1)^m\bigr(q_{s,m}\bigr)^{n-2}\Bigr],\quad s=1,2,
\end{equation*}
\begin{equation*}
M=Q_2\sum_{k=0}^\infty \bigr(q_{1,2k}\bigr)^{n-2}
+Q_1\sum_{k=0}^\infty\bigr(q_{2,2k+1}\bigr)^{n-2}.
\end{equation*}
Here, $q_{s,m}$'s and ${\mathbf{c}}_{s,m}$'s are defined by \eqnref{c1m}, \eqnref{q1m}, \eqnref{c2m}, \eqnref{q2m}.
\end{lem}

\pf
Let $$h_*(\mathbf{x})=Q_2\sum_{m=0}^\infty \frac{(-1)^m\bigr(q_{1,m}\bigr)^{n-2} }{|\mathbf{x}-\mathbf{c}_{1,m}|^{n-2}}
-Q_1\sum_{m=0}^\infty \frac{(-1)^m\bigr(q_{2,m}\bigr)^{n-2}}{|\mathbf{x}-\mathbf{c}_{2,m}|^{n-2}},$$
then $h_*$ satisfies \eqnref{eq:101} except the last condition, boundary integral conditions.

Note that
\begin{equation}\label{eqn:charge}
\frac{1}{(2-n)\omega_n}\int_{\p D_i}\pd{}{\nu}\frac{1}{|\mathbf{x}-\mathbf{p}|^{n-2}}d\sigma(\mathbf{x})
=\left\{
\begin{array}{ll}
0,\quad&\mbox{for }\mathbf{p}\in\RR^n\setminus\bar{D}_i,\\
1,&\mbox{for }\mathbf{p}\in D_i,
\end{array}\right.
\end{equation}
and we have
\begin{align*}
&\frac{1}{(2-n)\omega_n}\int_{\p D_1}\p_\nu h_*(\mathbf{x})d\sigma(\mathbf{x})\\
&=Q_2\sum_{k=0}^\infty \bigr(q_{1,2k}\bigr)^{n-2}+Q_1\sum_{k=0}^\infty\bigr(q_{2,2k+1}\bigr)^{n-2}\\
&=\sum_{k=0}^\infty  \bigr(q_{1,2k}\bigr)^{n-2}\sum_{k=0}^\infty  \bigr(q_{2,2k}\bigr)^{n-2}
-\sum_{k=0}^\infty\bigr(q_{1,2k+1}\bigr)^{n-2}\sum_{k=0}^\infty\bigr(q_{2,2k+1}\bigr)^{n-2}\\
&=Q_1\sum_{k=0}^\infty\bigr(q_{2,2k}\bigr)^{n-2}+Q_2\sum_{k=1}^\infty \bigr(q_{1,2k+1}\bigr)^{n-2}\\
&=-\frac{1}{(2-n)\omega_n}\int_{\p D_2}\p_\nu h_*(\mathbf{x})d\sigma(\mathbf{x}).
\end{align*}
\qed


\begin{lem}\label{Q1Q2} Assume that the dimension $n$ is $3$ and  the distance $\epsilon$ is sufficiently small. Then, there is a positive constant $C$ independent of $r_1$, $r_2$, and
$\epsilon$ satisfying the following properties:
\begin{itemize}
\item  {\textbf{Estimates for $ \sum_{m=0}^{\infty}q_{s,m} $:}}\begin{align}
\frac{1}{C}\frac d {d+1}|\log\epsilon|\label{sum:qtotal1} \ \leq\
&\sum_{m=0}^{\infty}q_{1,m}\ \leq \ C \frac d
{d+1}|\log\epsilon|,\\\label{sum:qtotal2} \frac{1}{C}\frac 1
{d+1}|\log\epsilon|\ \leq\ &\sum_{m=0}^{\infty}q_{2,m}\ \leq \
C\frac 1 {d+1}|\log\epsilon|
\end{align}
where $d=\frac{r_2}{r_1}.$
\item {\textbf{Estimates for $ Q_s $:}}

\begin{align*}
 \frac{1}{C} {\frac1 {d+1}} \  \leq Q_1  & \leq \  C{\frac 1 {d+1}},\\
 \frac{1}{C} {\frac d {d+1}}\  \leq Q_2  & \leq   \  C {\frac d {d+1}}
 \end{align*}
where $Q_1$ and $Q_2$ are defined in Lemma \ref{def:Q}.
\item  {\textbf{Estimates for $ \sum_{k=0}^\infty (c_{s,2k})\bigr(q_{s,2k}\bigr)-\sum_{k=0}^\infty (c_{s,2k+1})\bigr(q_{s,2k+1}\bigr) $:}}
$$\frac{r_s}{C}\leq  (-1)^{s+1}\left[\sum_{k=0}^\infty (c_{s,2k})\bigr(q_{s,2k}\bigr)-\sum_{k=0}^\infty (c_{s,2k+1})\bigr(q_{s,2k+1}\bigr)\right]\leq C r_s.$$

\end{itemize}
\end{lem}
Now, we are ready to prove Theorem \ref{thm2}.
\subsection {Proof of Theorem \ref{thm2}}
\noindent\textbf{Potential Difference}\\
 We first
consider the case of $H(x_1,x_2,x_3)=x_1$. Then, Lemma
\ref{series:h} implies
\begin{align}
\ds u\big|_{\partial D_1} - u\big|_{\partial D_2}
&=\frac{Q_2}{M}\sum_{k=0}^\infty (c_{1,2k})\bigr(q_{1,2k}\bigr)
+\frac{Q_1}{M}\sum_{k=0}^\infty (c_{2,2k+1})\bigr(q_{2,2k+1}\bigr)\notag\\
&\quad-\frac{Q_2}{M}\sum_{k=0}^\infty(
c_{1,2k+1})\bigr(q_{1,2k+1}\bigr) -\frac{Q_1}{M}\sum_{k=0}^\infty
(c_{2,2k})\bigr(q_{2,2k}\bigr), \label{def_u}
\end{align}
where  $c_{s,j}$ is the $x_1$-coordinate of $\mathbf{c}_{s,j}$ for
$\ s=1,2, j\in\NN$, i.e.,
\begin{equation} (c_{s,j},0,0)=\mathbf{c}_{s,j}.\notag
\end{equation}
Lemma \ref{Q1Q2} allows one to estimate four positive valued terms
in the right hand side of (\ref{def_u}) so that
\begin{align*}
M&=Q_2\sum_{k=0}^\infty q_{1,2k}
+Q_1\sum_{k=0}^\infty q_{2,2k+1}\\
&\simeq  \left[\frac {d^2} {(d+1)^2}+{\frac 1 {(d+1)^2}}\right]\bigr|\log\epsilon\bigr|\\
&\simeq |\log\epsilon|.
\end{align*}
In total, we obtain
\begin{align}\left| u\big|_{\partial D_1} - u\big|_{\partial D_2}\right|
&\simeq  \frac  {r_1 r_2}{r_1+r_2}\frac {1}{|\log
\epsilon|}.\nonumber
\end{align}

{In the case of $H=a_2 x_2 + a_3 x_3$, the integration
\eqnref{eq:102} is zero; all point charges of $h$ lie on $x_1$ axis.
Therefore, there is no potential difference between inclusions, and
we have established the estimate for the potential difference
between $D_1$ and $D_2$.}

Therefore, for $H= \sum_{i=1}^3 a_i x_i$,
\begin{align}\left| u\big|_{\partial D_1} - u\big|_{\partial D_2}\right|
&\simeq   |a_1|\frac {r_1 r_2}{r_1+r_2}\frac {1}{|\log
\epsilon|}.\nonumber
\end{align}

\noindent\textbf{Lower bound}\\
The lower bound is obtained { by simply applying} the mean value
theorem. Since $\left|u|_{D_1} - u|_{D_2}\right|$ behaves as
$1/|\log\epsilon|$, the gradient behaves as
$1/(\epsilon|\log\epsilon|)$, and, more precisely, there is a point
{$\mathbf{x}_0$} between $D_1$ and $D_2$ satisfying that
$$ \frac 1 C {|a_1|} \left(\frac {r_1 r_2}{r_1 + r_2} \right) {\frac 1
{ | \epsilon \log \epsilon |} }  \leq |{\nabla u (\mathbf{x}_0)}|$$
where $H= \sum_{i=1}^3 a_i x_i$.

\smallskip

\noindent\textbf{Upper bound}\\
The upper bound of the gradient is derived by applying the methods
presented by Bao et al \cite {BLY}.

We assume that $r_1 \geq r_2$, and let $$ \frac 1 {r_1} D_i = \left
\{ \mathbf{x} \in \mathbb{R}^3 \big|~ r_1\mathbf{x} \in D_i
\right\},\quad  i=1,~2.$$ Note that $\frac 1 {r_1} D_1$ is a unit
sphere. Define a bounded domain $\Omega$, containing
$\frac{1}{r_1}D_1$ and $\frac{1}{r_1}D_2$ independently of
$\epsilon$,  as the sphere $B_{4} (0,0,0)$.
 Now, for a solution $u$ to \eqref{eq:001} for $H(\mathbf{x}) = \sum_{i=1}^3 a_i x_i$, define the scaled function $\tilde u$ as
$$\widetilde{u}(\mathbf{x}):= \frac 1 {r_1} u (r_1 \mathbf{x}).$$
Then  $\widetilde{u}$ is
also the solution to \eqref{eq:001} for $H(\mathbf{x}) =
\sum_{i=1}^3 a_i x_i$ with $\frac 1 {r_1}D_1$ and $\frac 1 {r_1} D_2$ instead of $D_1$ and $D_2$.
The estimate for the difference of $u$ between $\partial D_1$ and
$\partial D_2$ in this theorem yields
\begin{align}\left| \widetilde{u}\big|_{\partial (\frac 1 {r_1} D_1)} - \widetilde{u}\big|_{\partial (\frac 1 {r_1} D_2)}\right|
&\simeq  |a_1|\frac {d} {1+d}\frac {1}{|\log \delta|},\nonumber
\end{align}
where $$d=\frac {r_2} {r_1}~ \mbox{and} ~\delta = \frac {\epsilon}
{r_1}.$$
By the maximum principle, we have
\begin{align} \norm {\widetilde{u}-H}_{L^{\infty}(\partial \Omega) } &\leq
\left| \widetilde{u}|_{\partial \frac 1 {r_1}D_1} -
\widetilde{u}|_{\partial \frac 1 {r_1} D_2} \right|
+2 \norm {H}_{L^\infty (\Omega)}\notag\\
& \leq C \left(  |a_1|\left(\frac {d}{1 + d} \right) {\frac 1 {|
\log \delta |} + |a_1 |+|a_2| + |a_3|} \right),\notag\end{align} and,
as a result, $$\norm {\widetilde{u}}_{L^{\infty}(\partial \Omega) }
\leq C' \left( |a_1|\left(\frac {d}{1+d} \right) {\frac 1 {| \log
\delta |} + |a_1| + |a_2| + |a_3|} \right).$$

To estimate $|\nabla \tilde u|$ on $\partial (\frac 1 {r_1} D_1 \cup
\frac 1 {r_1} D_2)$, we define $v_3$ as in \cite{BLY}:
\begin{equation}\notag
\quad \left\{
\begin{array}{ll}
\Delta v_3 = 0 ~&\mbox{in}~\Omega\setminus \overline{ (\frac 1 {r_1}D_1\cup \frac 1 {r_1} D_2)}\notag\\
v_3=0 ~&\mbox{on}~ \partial (\frac 1 {r_1} D_1\cup \frac 1 {r_1} D_2)\notag\\
v_3=-\widetilde{u} ~&\mbox{on}~ \partial \Omega\notag
\end{array}
\right.
\end{equation}
We draw the attention of readers to Lemma \ref{lem:Li}, \ref{lem:Y} which are modified from \cite{BLY} to fit our problem.
For a reader's convenient, we provide the proofs at the end of this section.
\begin{lem}(\cite{BLY}) \label{lem:Li}
There is a constant $C$ independent of $d$ and  $\epsilon$
such that
\begin{equation}\label{eq:a}\norm{\nabla(\widetilde{ u}+v_3 )}_{L^{\infty} (\Omega\setminus
(\frac 1 {r_1} D_1 \cup \frac 1 {r_1} D_2))} \leq C |a_1|\left(\frac
{d}{1+d} \right) {\frac 1 {|\delta \log \delta |}}. \end{equation}
\end{lem}

Now, in estimating $|\nabla \tilde u|$ on $\partial (\frac 1 {r_1} D_1 \cup
\frac 1 {r_1} D_2)$, it is remained to be derived an upper bound of $|\nabla v_3|$ on $\partial (\frac 1 {r_1} D_1 \cup
\frac 1 {r_1} D_2)$. To do that, we define the harmonic function $\rho$ in $\Omega \setminus (\frac 1 {r_1} D_1 \cup \frac 1 {r_1} D_2  )$ as in
\cite{BLY}:
\begin{equation}\notag
\quad  \left\{
\begin{array}{ll}
\ds\triangle \rho = 0~&~\mbox{in}~\Omega \setminus (\frac 1 {r_1} D_1 \cup \frac 1 {r_1} D_2  )\notag,\\
\ds\rho = 0 ~&~\mbox{on}~\partial (\frac 1 {r_1} D_1 \cup \frac 1 {r_1}
D_2  )\notag,\\
\ds\rho=1 &~\mbox{on}~\partial \Omega \notag.
\end{array}
\right.
\end{equation}
Note that $v_3=\pm\norm{\widetilde{u}}_{L^{\infty}(\partial \Omega)} \rho =
 0$ on $\partial \Bigr(\frac 1 {r_1} D_1 \cup \frac 1 {r_1} D_2\Bigr)$. Moreover, from the fact that $v_3 = - \widetilde{u}$ on the $\partial
\Omega$ and the maximum principle, for $x\in\Omega \setminus \Bigr(\frac 1 {r_1} D_1 \cup \frac 1 {r_1}
D_2\Bigr)$, we have $- \norm{\widetilde{u}}_{L^{\infty}(\partial
\Omega)} \rho\leq v_3 \leq
\norm{\widetilde{u}}_{L^{\infty}(\partial \Omega)}\rho$.
Therefore, by Hopf's Lemma and the maximum principle,
$$\norm{\nabla v_3 }_{L^{\infty} (\Omega\setminus
(\frac 1 {r_1} D_1 \cup \frac 1 {r_1} D_2))} \leq
\norm{\widetilde{u}}_{L^{\infty}(\partial \Omega)}  \norm{\nabla
\rho }_{L^{\infty} (\Omega\setminus (\frac 1 {r_1} D_1 \cup \frac 1
{r_1} D_2))}.$$
We apply the following lemma to calculate $\norm{\nabla v_3 }_{L^{\infty} (\Omega\setminus
(\frac 1 {r_1} D_1 \cup \frac 1 {r_1} D_2))}$.
\begin{lem}(\cite{BLY}) \label{lem:Y}
There is a constant $C$ such that $$\norm {\nabla \rho
}_{L^{\infty}(\partial \frac 1 {r_1} D_1 \cup \partial \frac 1 {r_1}
 D_2)} \leq C \frac 1 d,
$$ for $\epsilon$ small enough.
\end{lem}

 Applying Lemma \ref{lem:Y}, we have \be\norm {\nabla
v_3}_{L^{\infty}(\partial (\frac 1 {r_1} D_1 \cup \frac 1 {r_1}
D_2)) } \leq C \frac 1 d \left( |a_1|\left(\frac {d}{1+ d } \right)
{\frac 1 {| \log \delta |} +|a_1| +|a_2| + |a_3|}
\right).\label{eq:b}\ee Two bounds \eqref{eq:a} and \eqref{eq:b}
yield
\begin{align}
\norm {\nabla u}_{L^{\infty}(\partial (D_1 \cup D_2)) } &= \norm
{\nabla \widetilde{u}}_{L^{\infty}(\partial (\frac 1 {r_1}D_1 \cup
\frac 1 {r_1} D_2)) }\notag\\
&\leq C \left( |a_1|\left(\frac {d}{1 + d} \right) {\frac 1 {|
\delta \log \delta |} + \frac 1 d (|a_2| + |a_3|)} \right)\notag\\
&\leq C' \left( |a_1|\left(\frac {r_1 r_2}{r_1 + r_2} \right) {\frac
1 {| \epsilon \log \epsilon |} + \frac 1 d (|a_2| + |a_3|)}
\right)\notag\end{align} for sufficiently small $\epsilon>0$.
Since $|\nabla H|$ is bounded by $|a_1| + |a_2| + |a_3|$, we have
$$\norm {\nabla (u - H)}_{L^{\infty}(\partial (D_1 \cup D_2)) } \leq
C \left( |a_1|\left(\frac {r_1 r_2}{r_1 + r_2} \right) {\frac 1 {|
\epsilon \log \epsilon |} + \frac 1 d (|a_2| + |a_3|} )\right).$$ By
the harmonicity of $ u - H$ in $\mathbb{R}^3 \setminus (D_1 \cup
D_2)$, $\norm {\nabla (u - H)}_{L^{\infty}( \mathbb{R}^3 \setminus
(D_1 \cup D_2)) }$ has the same upper bound as the above. The fact
of $|\nabla H| \leq |a_1| + |a_2| + |a_3|$ is again used so that
$$\norm {\nabla u}_{L^{\infty}( \mathbb{R}^3 \setminus (D_1 \cup D_2)) }
\leq C \left( |a_1|\left(\frac {r_1 r_2}{r_1 + r_2} \right) {\frac 1
{| \epsilon \log \epsilon |} + \frac 1 d(|a_2| + |a_3|)} \right).$$
Therefore, we obtain the upper bound of the gradient estimate. \qed

\smallskip

\noindent\textbf{Proof of Lemma \ref{lem:Li}}\\
From definition, $\widetilde{u} + v_3$ is constant on $\partial
\frac 1 {r_1} D_1$ and $\frac 1 {r_1} D_1$ is a unit ball. By Kelvin
transform, $\widetilde{u} + v_3$ can be extended harmonically to
$\Omega \setminus (\frac 1 {r_1} D_{1,\delta'} \cup \frac 1 {r_1}
D_{2} )$ where $$ \frac 1 {r_1} D_{1,\delta'}  = \left \{\mathbf{ x}
\in \frac 1 {r_1} D_{1}\big| \mbox{dist}(\mathbf{x},
\partial {\frac 1 {r_1} D_1 } ) > \delta'\right\}
$$ and
$$\delta'=1 - \frac {1} {1+2\delta }.$$
Similarly, $\widetilde{u} + v_3$ can be also extended harmonically
to $\Omega \setminus (\frac 1 {r_1} D_{1,\delta'} \cup \frac 1 {r_1}
D_{2, \delta ''} )$ where $$ \frac 1 {r_1} D_{2,\delta''}  = \left
\{\mathbf{ x} \in \frac 1 {r_1} D_{2}\big| \mbox{dist}(\mathbf{x},
\partial {\frac 1 {r_1} D_2 } ) > \delta''\right\}
$$ and
$$\delta''=d - \frac {d^2} {d+2\delta }.$$
Furthermore, by the standard estimate for the extension, there is a
constant $C$ such that
$$\max_{\Omega \setminus (\frac 1 {r_1} D_{1,\delta'} \cup \frac 1 {r_1}
D_{2, \delta ''} ) } (\widetilde{u} - v_3) -\min_{\Omega \setminus
(\frac 1 {r_1} D_{1,\delta'} \cup \frac 1 {r_1} D_{2, \delta ''} ) }
(\widetilde{u} - v_3)  \leq C |a_1|\left(\frac {d}{1+d} \right)
{\frac 1 {|\log \delta |}}$$ for sufficiently small $\epsilon>0$.
Note that
$$\delta'\approx 2\delta~\mbox{and}~\delta''\approx 2\delta.$$
By the gradient estimate for harmonic functions, we have
$$\norm{\nabla(\widetilde{ u}+v_3 )}_{L^{\infty} (\Omega\setminus
(\frac 1 {r_1} D_1 \cup \frac 1 {r_1} D_2))} \leq C |a_1|\left(\frac
{d}{1+d} \right) {\frac 1 {|\delta \log \delta |}}. $$ \qed

\noindent\textbf{Proof of Lemma \ref{lem:Y}}\\
 Let $\rho_i$ $(i=1,2)$ be the solution to
\begin{equation}\notag
\quad  \left\{
\begin{array}{ll}
\triangle \rho_i = 0~&~\mbox{in}~\Omega \setminus (\frac 1 {r_1} D_i )\notag\\
\rho_i = 0 ~&~\mbox{on}~\partial (\frac 1 {r_1} D_i)\notag\\
\rho_i=1 &~\mbox{on}~\partial \Omega \notag
\end{array}
\right.
\end{equation}
Then $\rho = \rho_i$ on the $\partial \Omega \cup \partial D_i$. The
maximum principle yields to $\rho_{i} \leq \rho$. Note that the
radii of $\frac 1 {r_1} D_1$ and $\frac 1 {r_1} D_2$ are $1$ and $d$
respectively.

Consider the harmonic function $v$ which is the solution to
\begin{equation}\notag
\quad \left\{
\begin{array}{ll}
\Delta v = 0 ~&\mbox{in}~B_{4} (0,\cdots,0)\setminus B_{r_0} (\mathbf{c}_0)\notag,\\
v=0 ~&\mbox{on}~ \partial B_{r_0} (\mathbf{c}_0)\notag,\\
v=1 ~&\mbox{on}~ \partial B_4 (0,\cdots,0)\notag,
\end{array}
\right.
\end{equation}
where $r_0 \leq 1 $ and $|\mathbf{c}_0| \leq 2$. Let  the harmonic
function $w$ be as
\begin{equation}\notag
\quad w= \left\{
\begin{array}{ll}
\left(  \frac 1 {2^{n-2}} -\frac 1 {{r_0}^{n-2}}  \right)^{-1}
\left(  \frac 1 {|\mathbf{x}- \mathbf{c}_0|^{n-2}} -\frac 1 {{r_0}^{n-2}}  \right)~&~\mbox{for}~n \geq 3\notag\\
\left(\log 2 - \log {r_0}\right) \left(\log
{|\mathbf{x}-\mathbf{c}_0|} - \log
{r_0}\right)~&~\mbox{for}~n=2\notag
\end{array}
\right.
\end{equation}
Then, $$w \geq 1 =v~\mbox{on}~\partial B_{4} (0,\cdots,0)$$ and
$$w=v~\mbox{on}~\partial B_{r_0} (\mathbf{c}_0).$$
Then,
 there is a constant $C$, independent of $r_0$ and
$\mathbf{c}_0$, satisfying
$$ \norm { \nabla v}_{L^{\infty} (\partial B_{r_0} (\mathbf{c}_0))} \leq C \frac 1 {r_0}.$$
By Hopf's Lemma, we have
$$-\partial_{\nu} w \geq -\partial_{\nu}v \geq 0~\mbox{on}~{\partial B_{r_0} (\mathbf{c}_0)},$$
and there is a constant $C$, independent of $\mathbf{c}_0$ and
$r_0$, satisfying
$$-\partial_{\nu} w \leq C  \frac 1 {r_0}.$$
Therefore, we obtain
\begin{equation}\label{rho:1om} \norm { \nabla v}_{L^{\infty} (\partial
B_{r_0} (\mathbf{c}_0))}= -\partial_{\nu} w  \leq C \frac 1
{r_0}.\end{equation}

 It follows from \eqref{rho:1om} and Hopf's
Lemma that
$$\norm {\nabla \rho
}_{L^{\infty}(\partial \frac 1 {r_1} D_1)} \leq \norm {\nabla \rho_1
}_{L^{\infty}(\partial \frac 1 {r_1} D_1 )} \leq C$$ and
$$\norm {\nabla \rho
}_{L^{\infty}(\partial \frac 1 {r_1} D_2)} \leq \norm {\nabla \rho_2
}_{L^{\infty}(\partial \frac 1 {r_1} D_2 )} \leq C \frac 1 d.$$ \qed

\subsection{Proof of Lemma \ref{Q1Q2}}

We begin by considering the position sequences $\mathbf{c}_{1,m}$
and $\mathbf{c}_{2,m}$, because the quantities like $q_{i,m}$ and
$\rho_{i,j}$ used in the derivation are yielded by these position
sequences $\mathbf{c}_{i,m}$. Referring to the relation (\ref{c1m}),
the successor $\mathbf{c}_{1,m+2}$ to $\mathbf{c}_{1,m}$ is
determined by the twice refection $R_1 \circ R_2$. For
$\mathbf{x}=(x,0,\dots,0)\in D_1$, the twice reflected point
$R_1(R_2(\mathbf{x}))=(x',0,\dots,0)$ is given by
\begin{equation*}\label{eqn:R1R2}
x'=r_1+\epsilon-{\frac{r_1^2} {r_1+r_2+2\epsilon -{\frac {r_2^2}{x +
r_2+\epsilon}}}}.
\end{equation*}
For the sake of convenience, we assume that
 \be\label{def:y}
 \delta=\frac{\epsilon}{r_1},\ d=\frac{r_2}{r_1},\mbox{ and }y_{j}=\frac{c_{1,j}}{r_1},\quad
 j\in\NN,
 \ee
 where $c_{1,j}$ is the $x_1$-coordinate of $\mathbf{c}_{1,j}$.
Then, we have the  relation
\begin{align}
y_{2k}y_{2k-2}+\frac{d+(1+3d)\delta+2\delta^2}{1+d+2\delta}y_{2k}-&\frac{d+(3+d)\delta+2\delta^2}{1+d+2\delta}y_{2k-2}\nonumber\\
&\qquad-\frac{4d\delta+3(1+d)\delta^2+2\delta^3}{1+d+2\delta}=0.\label{yR1R2}
\end{align} Let $\mathbf{p}=(p,0,\dots,0)\in D_1$ be the fixed point of
$(R_1\circ R_2)$, then $\frac{p}{r_1}$ is the limit point of
$y_{2k}$ and  satisfies
\begin{equation*}\left(\frac{p}{r_1}\right)^2+\frac{2(d-1)\delta}{1+d+2\delta}\left(\frac{p}{r_1}\right)
-\frac{4d\delta+3(1+d)\delta^2+2\delta^3}{1+d+2\delta}=0,\end{equation*}
and, as a results,
\begin{equation}\label{papprox}
\frac{p}{r_1}=2\sqrt{\frac{d}{d+1}}\sqrt{\delta}+O(\delta).
\end{equation}
{Here, we have a constant $C$ independent of $d$ so that
$$|O(\sqrt \delta)| \leq C \sqrt \delta ~\mbox{for sufficiently small}~\delta>0.$$}

\par It follows from \eqnref{yR1R2} that
\begin{align*}
\Bigr(y_{2k}-\frac{p}{r_1}\Bigr)\Bigr(y_{2k-2}-\frac{p}{r_1}\Bigr)&
+\Bigr(\frac{d+(1+3d)\delta+2\delta^2}{1+d+2\delta}+\frac{p}{r_1}\Bigr)\Bigr(y_{2k}-\frac{p}{r_1}\Bigr)\\
&\qquad-\Bigr(\frac{d+{(3+d)}\delta+2\delta^2}{1+d+2\delta}-\frac{p}{r_1}\Bigr)\Bigr(y_{2k-2}-\frac{p}{r_1}\Bigr)=0.
\end{align*}
For simplicity, let
\begin{equation*}
z_{2k}=y_{2k}-\frac{p}{r_1},\end{equation*}
 then
\begin{equation*}\label{eqn:z}
1+\Big(\frac{d+(1+3d)\delta+2\delta^2}{1+d+2\delta}+\frac{p}{r_1}\Bigr)\frac{1}{z_{2k-2}}
-\Bigr(\frac{d+{(3+d)}\delta+2\delta^2}{1+d+2\delta}-\frac{p}{r_1}\Bigr)\frac{1}{z_{2k}}=0.
\end{equation*} Further, let $$ \ds A=
\frac{\left(\frac{d+(1+3d)\delta+2\delta^2}{1+d+2\delta}\right)+\frac{p}{r_1}}{\left(\frac{d+(3+d)\delta+2\delta^2}{1+d+2\delta}\right)-\frac{p}{r_1}}
\ \mbox{  and  }\
B=\frac{1}{2\left(\frac{(d-1)\delta}{1+d+2\delta}+\frac{p}{r_1}\right)},$$
this can be rewritten as
$$\Bigr(\frac{1}{z_{2k}}+B\Bigr)=A\Bigr(\frac{1}{z_{2k-2}}+B\Bigr).$$
Therefore, the sequence $y_{2k}$ is expressed as
\begin{align}
y_{2k}&=z_{2k}+\frac{p}{r_1}\notag\\&=\frac{1}{(\frac{1}{z_0}+B)A^k-B}+\frac{p}{r_1}\notag\\
&=\frac{1}{(\frac{1}{1+\delta-\frac{p}{r_1}}+B)A^k-B}+\frac{p}{r_1}.\label{complex_to_omuch}
\end{align}

\subsubsection{Estimates for $y_{j}$}
{{We simplify \eqref{complex_to_omuch}  under the assumption that
$\delta$ is sufficiently small. Our strategy is to choose an
appropriate $N$ so that $\sum_{k\leq N}$ is dominant in the
following series calculation, for example in
$\sum(q_{1,2m}-q_{1,2m+1})$, and estimate the series separately for
smaller and larger sub-indices.

From now on, we use the big $O$ notation frequently. The equation $f=g+O(\delta)$ means that there exist a constant $C$ independent of $\delta$ such that
$|f-g|\leq C\delta$ for small enough $\delta>0$. In this paper, $C$ is assumed additionally to be independent of $d$ and $k$ as well. We define $O(\sqrt{\delta})$ similarly.

 }}

The expression (\ref{complex_to_omuch}) of $y_{j}$ is too
complicated to well describe the dependency of $y_{j}$ with respect
to $d$ and $\delta$. Thus, a simplified expressions is established,
provided that the distance $\delta$ is small enough. As for our
strategy, we choose an appropriate number $N \simeq \sqrt \epsilon$
so that the sequence terms of $k\leq N$ are dominant in the sequence
$y_k$, and thus estimate $y_{2k}$ in two cases of $k\leq N$ and
$k\geq N$ separately.

\par From the definition of $A$ and $B$, we have
\begin{align*}
\ds A&=\frac{\frac{d}{1+d}+O(\delta)+\frac{p}{r_1}}{\frac{d}{1+d}+O(\delta)-\frac{p}{r_1}}=1+2\frac{d+1}{d}\frac{p}{r_1}+\frac{d+1}{d}O(\delta)\\
&=1+4\sqrt{\frac{d+1}{d}}\sqrt\delta+\frac{d+1}{d}O(\delta),
\end{align*}
and
\begin{align*}
\sqrt{\delta} B&=
\frac{\sqrt{\delta}}{4\sqrt{\frac{d}{d+1}}\sqrt{\delta}+O(\delta)}\\
&=\frac{1}{4\sqrt{\frac{d}{d+1}}}+\frac{d+1}{d}O(\sqrt{\delta}).
\end{align*}
By a standard argument, one can show that for $x\in (0,2)$ and
$(1+x)^k\leq2$,
$$(1+x)^k\leq 1+kx+k^2x^2.$$
For $k\leq\frac{\log
2}{8}\sqrt{\frac{d}{d+1}}\frac{1}{\sqrt{\delta}}$, we have
$\Bigr(1+8\sqrt{\frac{d+1}{d}}\sqrt\delta\Bigr)^k\leq 2.$ This
yields
\begin{align*}A^k
&=1+k\Bigr(4\sqrt{\frac{d+1}{d}}\sqrt\delta+\frac{d+1}{d}O(\delta)\Bigr)+\frac{d+1}{d}k^2O(\delta)\\
&=1+4k\sqrt{\frac{d+1}{d}}\sqrt\delta+\frac{d+1}{d}k^2O(\delta).
\end{align*}
Hence, the estimate implies
\begin{align} \ds
y_{2k}&\ds=\frac{1}{\left(B+1+O(\sqrt\delta)\right)\left(1+4k\sqrt{\frac{d+1}{d}}\sqrt\delta+\frac{d+1}{d}k^2O(\delta)\right)-B}+\frac{p}{r_1}\nonumber\\\nonumber
\ds&=\frac{1}{1+k\frac{d+1}{d}+\bigr(\frac{d+1}{d}\bigr)^{\frac{3}{2}}k^2
O_1(\sqrt\delta)}+\frac{p}{r_1}.
\end{align}
Here is a constant $C_1 > 0$ independent of $\delta$, $d$ and $k$
such that
$$|O_1(\sqrt\delta)| \leq C_1 \sqrt \delta .$$

\par We take the integer $N$ as
$$N = \min_{n\in \mathbb{N}} \left\{n \geq \frac 1 {C_1 } \frac {\log 2 }{8 } \sqrt {\frac
{d}{d+1}} \frac {1}{\sqrt \delta}\right\}.$$ Then, for $k \leq N$,
\begin{align}
\ds
y_{2k}&\ds=\frac{d}{k(d+1)+d}+\sqrt{\frac{d}{d+1}}O(\sqrt\delta),\label{y2k}
\end{align}
and
\begin{align}
y_{2k+1}&=\frac{1}{r_1}\Bigr(-r_2-\epsilon+\frac{r_2^2}{r_1y_{2k}+r_2+\epsilon}\Bigr)\nonumber\\\nonumber
&=-d-\delta+\frac{d^2}{\frac{d}{k(d+1)+d}+d+\sqrt{\frac{d}{d+1}}O(\sqrt\delta)}\\
&=-\frac{d}{d+1}\frac{1}{k+1}+\sqrt{\frac{d}{d+1}}O(\sqrt\delta).\label{y2k1}
\end{align}

\par In the case of $k\geq N$, we use the fact
that the sequence $y_{1,2k}$ is decreasing to $\frac{p}{r_1} $,
i.e., $y_{0}>y_{2}\ \
>\cdots\ >\ y_{2k}>y_{2(k+1)}>\cdots>\frac{p}{r_1}.$
Here \eqnref{y2k} yields $y_{1,2N} =
\sqrt{\frac{d}{d+1}}O(\sqrt\delta) $. From the estimate
\eqnref{papprox} for $\frac p {r_1}$, we have
\begin{equation}\label{yp:approx1}
C\sqrt{\frac{d}{d+1}}\sqrt{\delta}>y_{2k}>\sqrt{\frac{d}{d+1}}\sqrt{\delta}
\end{equation}
and
\begin{equation}\label{yp:approx2}
C\sqrt{\frac{d}{d+1}}\sqrt{\delta}>-y_{2k+1}>\sqrt{\frac{d}{d+1}}\sqrt{\delta}.\end{equation}
Therefore, all estimates (\ref {y2k}), (\ref{y2k1}), (\ref
{yp:approx1}) and (\ref {yp:approx2}) for $y_{2k} $ and $y_{2k+1}$
are obtained.

\subsubsection{Estimates for $ \sum_{m=0}^{\infty}q_{s,m} $}

\par We consider the estimate for $$q_{1,m}=\prod_{j=0}^{m}\rho_{1,j}.$$
For $k \leq N$,  from \eqnref{rho1}, \eqnref{y2k} and \eqnref{y2k1},
we have
\begin{align}
\rho_{1,2k+1}&=\frac{d}{d+\delta+y_{2k}}=\frac{k(d+1)+d}{(k+1)(d+1)}\bigr(1+\frac{1}{\sqrt{d(d+1)}}O(\sqrt{\delta}\bigr)),\label{rho12}\\
\rho_{1,2k+2}&=\frac{1}{1+\delta-y_{2k+1}}=\frac{(k+1)(d+1)}{(k+1)(d+1)+d}\bigr(1+\sqrt{\frac{d}{d+1}}O(\sqrt{\delta})\bigr),\nonumber
\end{align}
and$$(\rho_{1,2k+1})(\rho_{1,2k+2})=\frac{k(d+1)+d}{(k+1)(d+1)+d}\bigr(1+\sqrt{\frac{d+1}{d}}O(\sqrt{\delta})\bigr).$$
In the case of $m \leq N$, they lead to
\begin{align}
q_{1,2m}&=\Bigr(\prod_{k=0}^{m-1}\rho_{1,2k+1}\rho_{1,2k+2}\Bigr)=\frac{d}{m(d+1)+d}\Bigr(1+\sqrt{\frac{d+1}{d}}O\bigr(\sqrt{\delta} \bigr)\Bigr)^m\label{approx:q1}\\\nonumber
&=\frac{d}{m(d+1)+d}\Bigr(1+m\sqrt{\frac{d+1}{d}}O\bigr(\sqrt{\delta}\bigr)\Bigr),\\\nonumber
q_{1,2m+1}&=(q_{1,2m})(\rho_{1,2m+1})\\
&=\frac{d}{(m+1)(d+1)}\Bigr(1+m\sqrt{\frac{d+1}{d}}O\bigr(\sqrt{\delta}\bigr)\Bigr).\nonumber
\end{align}
Thus, there exists a positive constant $C$ such that, for $s=0,1$,
\begin{align*}\frac{1}{C}\frac{d}{d+1}\frac{1}{m+1}\ \leq \ & q_{1,2m+s}\ \leq\  C\frac{d}{d+1}\frac{1}{m+1},\quad 1\leq m\leq N
\end{align*}
and, therefore,
\begin{align*}\label{sum:qsmall}
\frac{1}{\tilde{C}}\frac{d}{d+1}|\log\delta|\ \leq\ & \sum_{m\leq N}\bigr(q_{1,2m}+q_{1,2m+1}\bigr)
\ \leq\  \tilde{C}\frac{d}{d+1}|\log\delta|.
\end{align*}

\par In the case of $m\geq N$,  \eqnref{yp:approx1} and
\eqnref{yp:approx2}  yield, for $s=1,2$,
\begin{equation*}\label{property:rho}
\frac{1}{C}\frac{1}{1+\sqrt{\frac{\min(1/d,d)}{1+d}}\sqrt\delta}\ \leq\ \rho_{s,m}\leq C\frac{1}{1+\sqrt{\frac{\min(1/d,d)}{1+d}}\sqrt\delta},\quad m\geq N.
\end{equation*} Then it leads to
\begin{align}
\ds\sum_{m\geq N}\bigr(q_{1,2m}+q_{1,2m+1}\bigr)\nonumber
\ds&\simeq\bigr(q_{1,2N}+q_{1,2N+1}\bigr)\sum_{j=1}^\infty \Bigr(\frac{1}{1+\sqrt{\frac{\min(1/d,d)}{1+d}}\sqrt\delta}\Bigr)^j\\\nonumber
\ds&\simeq \bigr(q_{1,2N}+q_{1,2N+1}\bigr)\frac{1}{\sqrt{\frac{\min(1/d,d)}{1+d}}\sqrt\delta}\nonumber\\\nonumber
&\simeq  \frac{d}{d+1}\frac{1}{N}\frac{1}{\sqrt{\frac{\min(1/d,d)}{1+d}}\sqrt\delta}\\
\ds&\simeq  \max(d,1).\label{sum:qlarge}
\end{align}
Therefore, we obtain \eqnref{sum:qtotal1}, and replacing $d$ by
$\frac{1}{d}$ we also have \eqnref{sum:qtotal2}.
\smallskip

\subsubsection{Estimates for $Q_s$}
We consider $Q_1=\sum_{m=0} ^{\infty} \bigr(q_{1,2m}-q_{1,2m+1}\bigr
)$. From definition, $q_{1,m}$ has the decreasing property as
\begin{align*}
\ds q_{1,2m+1}&=(\rho_{1,2m+1})(q_{1,2m})< q_{1,2m},\\
\ds q_{1,2m}&=(\rho_{1,2m})(q_{1,2m-1})< q_{1,2m-1},
\end{align*}
 and therefore
\begin{align}\label{Q1large}
0<\sum_{m\geq {N}}\bigr(q_{1,2m}-q_{1,2m+1}  \bigr)&< q_{1,2{N}}\leq
C \frac{d}{(d+1)(N+1)}.
\end{align}
This means that $\sum_{m\geq {N}}\bigr(q_{1,2m}-q_{1,2m+1} \bigr)$
shrinks to $0$ as $\delta$ goes to $0$.

\par On the other hand, it follows from  (\ref{rho12})  and
(\ref{approx:q1}) that

\begin{align*}
\ds\sum_{0\leq m<N} \bigr(q_{1,2m}-q_{1,2m+1} \bigr) &= \sum_{m<N}
q_{1,2m} \left(1 - \rho_{1,2m+1} \right)\\ &= \Bigr[\sum_{m<N}
\frac{d}{(m+1)(d+1)\bigr(m(d+1)+d\bigr)}\Bigr]\Bigr(1+N\sqrt{\frac{d+1}{d}}O\bigr(\sqrt{\delta}\bigr)\Bigr).
\end{align*}
Taking an advantage of a strict decreasing sequence, we get
\begin{align*}
\frac{1}{2}\Bigr(\frac{1}{d+1}+\mathfrak{A}\Bigr)<\sum_{0\leq m<N}
&\frac{d}{(m+1)(d+1)\bigr(m(d+1)+d\bigr)}<
\frac{1}{d+1}+\mathfrak{A},
\end{align*}
where
\begin{align*}
\mathfrak{A}&=\int_0^{M}\frac{d}{(t+1)(d+1)\bigr(t(d+1)+d\bigr)}\ dt\\
&=\frac{d}{d+1} \int_0^{{N}}\left[\frac{1}{t+\frac{d}{d+1}}-\frac{1}{t+1}\right]\ dt\\
&=\frac{d}{d+1}\Bigr[\log\frac{d+1}{d}+\log\bigr(\frac{N+\frac{d}{d+1}}{N+1}\bigr)\Bigr].
\end{align*}
As has been mentioned, $\delta$ is assumed to be small enough so
that
$$0\leq -\log\bigr(\frac{N+\frac{d}{d+1}}{N+1}\bigr)\leq\frac{1}{2}\log\frac{d+1}{d}.$$
It leads to
 \begin{align*}
 \frac{1}{4} \frac{d}{d+1}
\bigr(\log\frac{d+1}{d}+ \frac 1 d \bigr)\  \leq Q_1 & \leq \frac{3}{2}  \frac{d}{d+1} \bigr(\log\frac{d+1}{d}+ {\frac 1 d} \bigr).
 \end{align*}
Note that $$\log(1+x)<x,\quad\mbox{for }x>0.$$ Therefore, we obtain
 \begin{align*}
 \frac{1}{C_1} \frac{1}{d+1}\  \leq Q_1 & \leq  C_1 \frac{1}{d+1} .
 \end{align*}
Similarly, replacing $d$ by $\frac 1 d$, we also have
 \begin{align*}
 \frac{1}{C_2} \frac{d}{d+1}\  \leq Q_2 & \leq  C_2 \frac{d}{d+1} .
 \end{align*}

\subsubsection{Estimates for  $\sum_{k=0}^\infty (c_{s,2k})\bigr(q_{s,2k}\bigr)-\sum_{k=0}^\infty (c_{s,2k+1})\bigr(q_{s,2k+1}\bigr)$}

 \par The lower and upper bounds of
$$\sum_{m=0}^\infty (c_{s,2m})\bigr(q_{s,2m}\bigr) -
\sum_{m=0}^\infty( c_{s,2m+1})\bigr(q_{s,2m+1}\bigr),~~ s=1,2$$ are
established here. From \eqnref{def:y}, \eqnref{y2k}, \eqnref{y2k1}
and \eqnref{sum:qlarge}, we calculate
\begin{align*}
\frac{1}{r_1}\ds\sum_{m=0}^\infty (c_{1,2m})\bigr(q_{1,2m}\bigr)
\ds&=\sum_{0\leq m\leq N} (y_{1,2m})\bigr(q_{1,2m}\bigr)+\sum_{m>N} (y_{1,2m})\bigr(q_{1,2m}\bigr)\\
\ds&=\sum_{0\leq m\leq N} C\Bigr(\frac{d}{m(d+1)+d}\Bigr)^2+\sqrt{\frac{d+1}{d}}O(\sqrt\delta)
\sum_{m>N} q_{1,2m}\\
\ds&\leq C.
\end{align*}
Moreover,
\begin{align*}
\frac{1}{r_1}\ds\sum_{k=0}^\infty (c_{1,2k})\bigr(q_{1,2k}\bigr)\geq (y_{1,0})(q_{1,0})=1.
\end{align*}
Similarly, we have
$$(y_{1,1}) (-q_{1,1})\leq -\frac{1}{r_1}\ds\sum_{k=0}^\infty (c_{1,2k+1})\bigr(q_{1,2k+1}\bigr)\leq C,$$
and, in total, we have
$$\frac{1}{C}\leq\frac{1}{r_1}\left[\sum_{k=0}^\infty (c_{1,2k})\bigr(q_{1,2k}\bigr)-\sum_{k=0}^\infty (c_{1,2k+1})\bigr(q_{1,2k+1}\bigr)\right]\leq C.$$
By the same way,
$$\frac{1}{C}\leq-\frac{1}{r_2}\left[\sum_{k=0}^\infty (c_{2,2k})\bigr(q_{2,2k}\bigr)-\sum_{k=0}^\infty (c_{2,2k+1})\bigr(q_{2,2k+1}\bigr)\right]\leq C.$$
 \qed

\subsection{The Derivation for Theorem \ref{thm3}}

\begin{lem}\label{nover4} Assume that the dimension $n \geq 4$ and  the distance $\epsilon$ is sufficiently small. Then, there is a positive constant $C$ independent of $r_1$, $r_2$, and
$\epsilon$ satisfying the following properties:
\begin{itemize}
\item  {\textbf{Estimates for $ \sum_{m=0}^{\infty}q_{s,m}^{n-2}
$:}}$$ \sum_{m=0}^{\infty}q_{1,m}^{n-2} \simeq
1~\mbox{and}~\sum_{m=0}^{\infty}q_{2,m}^{n-2} \simeq 1
$$
where $d=\frac{r_2}{r_1}.$
\item {\textbf{Estimates for $ Q_s $:}}
$$ Q_1 \simeq {\frac 1 {d+1}} ~\mbox{and}~ Q_2  \simeq   {\frac d {d+1}}$$
where $Q_1$ and $Q_2$ are defined in Lemma \ref{def:Q}.
\item  {\textbf{Estimates for  $\sum_{k=0}^\infty (c_{s,2k})\bigr(q_{s,2k}\bigr)-\sum_{k=0}^\infty (c_{s,2k+1})\bigr(q_{s,2k+1}\bigr)$:}}
$$(-1)^{s+1}\left[\sum_{k=0}^\infty
(c_{s,2k})\bigr(q_{s,2k}\bigr)^{n-2}-\sum_{k=0}^\infty
(c_{s,2k+1})\bigr(q_{s,2k+1}\bigr)^{n-2}\right] \simeq
r_s~~\mbox{for}~s= 1,2. $$
\end{itemize}
\end{lem}

\pf Let $N$ be as chosen in the proof of Lemma \ref{Q1Q2}. We have
shown in \textbf{Estimates for} $\sum_{m=0}^{\infty} q_{s,m} $ that
for $j=1,2$,
\begin{align*}
q_{1,2m+j} \simeq\frac{d}{d+1}\frac{1}{m+1}\quad \mbox{for}~ 1\leq
m\leq N
\end{align*}
and $q_{1,0} = 1$. This yields
$$\sum_{0\leq m \leq N} q_{s,m}^{n-2} \simeq 1 $$
and by the argument of (\ref {sum:qlarge}), we have
\begin{align}
\ds\sum_{m\geq
N}\bigr(q_{1,2m}^{n-2}+q_{1,2m+1}^{n-2}\bigr)\nonumber
\ds&\simeq\bigr(q_{1,2N}^{n-2}+q_{1,2N+1}^{n-2}\bigr)\sum_{j=1}^\infty
\Bigr(\frac{1}{1+\sqrt{\frac{\min(1/d,d)}{1+d}}\sqrt\delta}\Bigr)^{j(n-2)}\\\nonumber
\ds&\simeq \bigr(q_{1,2N}^{n-2}+q_{1,2N+1}^{n-2}\bigr)
\frac{1}{\sqrt{\frac{\min(1/d,d)}{1+d}}\sqrt\delta}\nonumber\\\nonumber
&\simeq \left( \frac{d}{d+1}\frac{1}{N} \right)^{n-2}
\frac{1}{\sqrt{\frac{\min(1/d,d)}{1+d}}\sqrt\delta}.
\end{align}
Since $N$ increase at the rate of $\frac 1 {\sqrt {\delta}}$ as
$\delta$ goes to zero, $\sum_{m\geq
N}\bigr(q_{1,2m}^{n-2}+q_{1,2m+1}^{n-2}\bigr)$ shrinks to zero as
$\delta$ goes to zero. Therefore, we obtain $
\sum_{m=0}^{\infty}q_{1,m}^{n-2} \simeq 1$, and similarly
$\sum_{m=0}^{\infty}q_{2,m}^{n-2} \simeq 1 $.

\par Now we consider $Q_s = \sum_{m=0}^{\infty} (-1)^m (q_{s,m})^{n-2}$ for
$s=1,2$. We note that the tail sum $$ |\sum_{m \geq 2 N} (-1)^m
(q_{1,m})^{n-2}  | \leq  |\sum_{m \geq N} q_{1,2m}^{n-2}+
q_{1,2m+1}^{n-2} |
$$ and the upper bound in the right hand side has been shown above
to shrink to zero as $\delta$ goes to zero. In this respect, we
estimate only $\sum_{m \leq 2 N - 1} (-1)^m (q_{1,m})^{n-2}$ for
$Q_1$. From definition, we have
\begin{align}
\sum_{m \leq 2 N-1} (-1)^m (q_{1,m})^{n-2} &= \sum_{m \leq  N}
(q_{1,2 m})^{n-2} \left( 1 - (\rho_{2,2m+1})^{n-2}\right)\notag\\
& \simeq \sum_{m \leq  N} \left( \frac d {m(d+1) + d}\right)^{n-2}
{\left( 1 -
\rho_{2,2m+1}\right)}\notag \\
& \simeq \sum_{m \leq  N} \left( \frac d {m(d+1) + d}\right)^{n-2}
\left( \frac 1 {(m+1)(d+1)} + \frac 1 {\sqrt{d (d+1)}} O(\sqrt
{\delta})\right)\notag  \\
& \simeq \sum_{m \leq  N} \left( \frac d {m(d+1) + d}\right)^{n-2}
\left( \frac 1 {(m+1)(d+1)} \right)\notag
 \\
& \simeq \frac 1 {d+1} + \sum_{1< m \leq N}\frac
{d^{n-2}}{(d+1)^{n-1}(m+1)^{n-2}}\notag\\ & \simeq \frac 1
{d+1}.\notag
\end{align}
Therefore, we obtain $Q_1 \simeq \frac {1}{d+1}$, and replacing $d$
by $\frac 1 d$, also have $Q_2 \simeq \frac {d}{d+1}$.

\par We consider the last estimate. By Lemma \ref{Q1Q2}, we have
\begin{align}
1 &\leq \frac {c_{1,0}}{r_1} q_{1,0}^{n-1} \notag\\&\leq
\frac{1}{r_1}\left[\sum_{k=0}^\infty
(c_{1,2k})\bigr(q_{1,2k}\bigr)^{n-2}-\sum_{k=0}^\infty
(c_{1,2k+1})\bigr(q_{1,2k+1}\bigr)^{n-2}\right]\notag \\&\leq
\frac{1}{r_1}\left[\sum_{k=0}^\infty
(c_{1,2k})\bigr(q_{1,2k}\bigr)-\sum_{k=0}^\infty
(c_{1,2k+1})\bigr(q_{1,2k+1}\bigr)\right]\simeq 1.\notag
\end{align}
Therefore, we have
$$\left[\sum_{k=0}^\infty
(c_{1,2k})\bigr(q_{1,2k}\bigr)^{n-2}-\sum_{k=0}^\infty
(c_{1,2k+1})\bigr(q_{1,2k+1}\bigr)^{n-2}\right] \simeq r_1$$ and
similarly
$$-\left[\sum_{k=0}^\infty
(c_{2,2k})\bigr(q_{2,2k}\bigr)^{n-2}-\sum_{k=0}^\infty
(c_{2,2k+1})\bigr(q_{2,2k+1}\bigr)^{n-2}\right] \simeq r_2$$ \qed
\subsubsection {The proof of Theorem \ref{thm3}}
 We first
consider the case of $H(x_1,x_2,\cdots,x_n)=x_1$. Then, Lemma
\ref{series:h} implies
\begin{align}
\ds u\big|_{\partial D_1} - u\big|_{\partial D_2}
&=\frac{Q_2}{M}\sum_{k=0}^\infty
(c_{1,2k})\bigr(q_{1,2k}\bigr)^{n-2}
+\frac{Q_1}{M}\sum_{k=0}^\infty (c_{2,2k+1})\bigr(q_{2,2k+1}\bigr)^{n-2}\notag\\
&\quad-\frac{Q_2}{M}\sum_{k=0}^\infty(
c_{1,2k+1})\bigr(q_{1,2k+1}\bigr)^{n-2}
-\frac{Q_1}{M}\sum_{k=0}^\infty
(c_{2,2k})\bigr(q_{2,2k}\bigr)^{n-2}. \label{def_u2}
\end{align}
where  $c_{s,j}$ is the $x_1$-coordinate of $\mathbf{c}_{s,j}$ for
$\ s=1,2, j\in\NN$, i.e.,
\begin{equation} (c_{s,j},0,\cdots,0)=\mathbf{c}_{s,j}.\notag
\end{equation}

\par Lemma \ref{nover4} allows one to estimate four positive valued terms in the right hand side of
(\ref{def_u2}) so that
\begin{align*}
M&=Q_2\sum_{k=0}^\infty q_{1,2k}
+Q_1\sum_{k=0}^\infty q_{2,2k+1}\\
&\simeq  1
\end{align*}
In total, we obtain
\begin{align}\left| u\big|_{\partial D_1} - u\big|_{\partial D_2}\right|
&\simeq {\frac{r_1 r_2} {r_1+r_2} }.\nonumber
\end{align}
In the case of $H=\sum_{i=2}^n a_i x_i$, the integration
\eqnref{eq:102} is zero, because all point charges of $h$ lie on
$x_1$ axis. Thus, there is no potential difference between
inclusions. Therefore, we have established the estimate for the
potential difference between $D_1$ and $D_2$.

\par Now, we consider the upper bound of $|\nabla u|$ when $H(\mathbf{x})= \sum_{i=1} ^n a_i x_i$. To do so,
we pursue the argument similar to Theorem \ref{thm2}. We thus assume
that $r_1 \geq r_2$. Let $ \frac 1 {r_1} D_1$ and  $ \frac 1 {r_1}
D_2$ be  the same as defined in the proof of Theorem \ref{thm2}. We
choose the sphere {$B_{4} (0,\dots,0)$} as the domain $\Omega$
containing $D_1$ and $D_2$ independently of any small distance
$\epsilon>0$. We consider
$$\widetilde{u}(\mathbf{x}):= \frac 1 {r_1} u (r_1 \mathbf{x}).$$
Then, $\widetilde{u}$ is the solution to \eqref{eq:001} for
$H(\mathbf{x}) = \sum_{i=1}^3 a_i x_i$ when the inclusions are
$\frac 1 {r_1}D_1$ and $\frac 1 {r_1} D_2$ instead of $D_1$ and
$D_2$.  It follows from the estimate for the difference of $u$
between $D_1$ and $D_2$ that
\begin{align}\left| \widetilde{u}\big|_{\partial (\frac 1 {r_1} D_1)} - \widetilde{u}\big|_{\partial (\frac 1 {r_1} D_2)}\right|
&\simeq  |a_1|{\frac {d} {1+d}},\nonumber
\end{align}
where $d=\frac {r_2} {r_1}~ \mbox{and} ~\delta = \frac {\epsilon}
{r_1}.$  As defined to (0.9) in \cite{BLY}, $v_3$ is defined as
follows:
\begin{equation}\notag
\quad \left\{
\begin{array}{ll}
\Delta v_3 = 0 ~&\mbox{in}~\Omega\setminus \overline{ (\frac 1 {r_1}D_1\cup \frac 1 {r_1} D_2)}\notag\\
v_3=0 ~&\mbox{on}~ \partial (\frac 1 {r_1} D_1\cup \frac 1 {r_1} D_2)\notag\\
v_3=-\widetilde{u} ~&\mbox{on}~ \partial \Omega\notag
\end{array}
\right.
\end{equation}
Lemma \ref{lem:Li} and \ref{lem:Y} have been used in the proof of
Theorem \ref{thm2}. By the help of them, we can obtain \be
\norm{\nabla(\widetilde{ u}+v_3 )}_{L^{\infty} (\Omega\setminus
(\frac 1 {r_1} D_1 \cup \frac 1 {r_1} D_2))} \leq C |a_1|\left(\frac
{d}{1+d} \right) {\frac 1 {|\delta  |}}. \notag\ee and
 \be\norm
{\nabla v_3}_{L^{\infty}(\partial (\frac 1 {r_1} D_1 \cup \frac 1
{r_1} D_2)) } \leq C \frac 1 d \left( |a_1|\left(\frac {d}{1+ d }
\right)
 +\sum_{i=1} ^{n} |a_i|  \right).\notag\ee Two bounds lead to
\begin{align}
\norm {\nabla u}_{L^{\infty}(\partial (D_1 \cup D_2)) } &= \norm
{\nabla \widetilde{u}}_{L^{\infty}(\partial (\frac 1 {r_1}D_1 \cup
\frac 1 {r_1} D_2)) }\notag\\
&\leq C \left( |a_1|\left(\frac {d}{1 + d} \right) {\frac 1 {|
\delta|} + \frac 1 d \sum_{i=2} ^{n}|a_i| } \right)\notag\\
&\leq C' \left( |a_1|\left(\frac {r_1 r_2}{r_1 + r_2} \right) {\frac
1 {| \epsilon |} + \frac 1 d \sum_{i=2} ^n |a_i|} \right)\notag
\end{align}
for any sufficiently small $\epsilon>0$. Since $|\nabla H|$ is
bounded by $\sum_{i=1} ^{n} |a_i| $, the harmonicity of $ u - H$ in
$\mathbb{R}^n \setminus \Omega$ leads to
$$\norm {\nabla u}_{L^{\infty}( \mathbb{R}^n \setminus (D_1 \cup D_2)) }
\leq C \left( |a_1|\left(\frac {r_1 r_2}{r_1 + r_2} \right) {\frac 1
{| \epsilon |} + \frac 1 d \sum_{i=2}^n |a_i| } \right).$$
Therefore, we obtain the upper bound of the gradient estimate. In
addition, the lower bound of the gradient is immediately derived
from the mean value theorem. \qed

\subsection{The proof for Theorem \ref{thm4}}

Without the loss of generality, we may assume that $H(0,\cdots,0)=0$. Then, Lemma
\ref{series:h} implies
\begin{align}
\ds \left| u\big|_{\partial D_1} - u\big|_{\partial D_2} \right| &=
\Big| \frac{Q_2}{M}\sum_{k=0}^\infty
H(\mathbf{c}_{1,2k})\bigr(q_{1,2k}\bigr)^{n-2}
+\frac{Q_1}{M}\sum_{k=0}^\infty H(\mathbf{c}_{2,2k+1})\bigr(q_{2,2k+1}\bigr)^{n-2}\notag\\
&~~\quad-\frac{Q_2}{M}\sum_{k=0}^\infty H(
\mathbf{c}_{1,2k+1})\bigr(q_{1,2k+1}\bigr)^{n-2}
-\frac{Q_1}{M}\sum_{k=0}^\infty
H(\mathbf{c}_{2,2k})\bigr(q_{2,2k}\bigr)^{n-2}\Big| \notag\\
&\leq \norm {\nabla H}_{L^{\infty}(\Omega)} \Big|
\frac{Q_2}{M}\sum_{k=0}^\infty (c_{1,2k})\bigr(q_{1,2k}\bigr)^{n-2}
+\frac{Q_1}{M}\sum_{k=0}^\infty (c_{2,2k+1})\bigr(q_{2,2k+1}\bigr)^{n-2}\notag\\
&~~~~~~~~~~~~~~~~~~~\quad-\frac{Q_2}{M}\sum_{k=0}^\infty(
c_{1,2k+1})\bigr(q_{1,2k+1}\bigr)^{n-2}
-\frac{Q_1}{M}\sum_{k=0}^\infty
(c_{2,2k})\bigr(q_{2,2k}\bigr)^{n-2}\Big|. \notag
\end{align}
Owing to Lemma \ref{Q1Q2} and \ref{nover4}, we can choose a constant
$C$ independent of $\epsilon$, $r_1$, $r_2$ and $\norm {\nabla
H}_{L^{\infty}(\Omega)}$ satisfying

\begin{equation}\notag
\left| u\big|_{\partial D_1} - u\big|_{\partial D_2} \right| \leq C
\norm {\nabla H}_{L^{\infty}(\Omega)} \left(\frac {r_1 r_2}{r_1 +
r_2} \right)  \cdot \left\{
\begin{array}{ll}
\frac 1 {|\log \epsilon|}~&\mbox{if}~ n= 3\notag\\
1~&\mbox{if}~ n \geq 4.\notag
\end{array}
\right. \end{equation}

\par Here, we would establish the upper bound of the gradient.
The method is the same as presented in the proofs of Theorem
\ref{thm2} and \ref{thm3}. We thus assume that $r_1 \geq r_2$. Let
$\frac 1 {r_1} D_i~(i=1,~2)$, $\widetilde{u}$, $d$ and $\delta$ be
as defined in Theorem \ref{thm2} and \ref{thm3}. Then
$\widetilde{u}$ is the solution to \eqref{eq:001} for
$$\widetilde{H}:=\frac 1 r_1  H (r_1 \mathbf{x})$$ instead of $H$,
replacing $D_i$ by $\frac 1 {r_1} D_i$ for $i=1,2$.  We consider
$\frac 1 {r_1} \Omega := \left\{ \mathbf{x} \in \mathbb{R}^n~\big|~
r_1\mathbf{x}\in \Omega\right\}$, because $\Omega$ is the ball
$B_{4(r_1+r_2)} (0,\cdots,0)$ in this theorem. So, the radius of the
ball $\frac 1 {r_1} \Omega$ is between $4$ and $8$. Note that
$$\norm{\nabla\widetilde{H}}_{L^{\infty} (\Omega)} = \norm{\nabla{H}}_{L^{\infty} (\frac 1 {r_1}\Omega)} . $$
We thus choose the harmonic function $v_3$ satisfying the same
definition as Theorem \ref{thm2} and \ref{thm3}, replacing $\Omega$
by $\frac 1 {r_1} \Omega$. Similarly to Theorem \ref{thm2} and
\ref{thm3}, owing to Lemma \ref{lem:Li}and \ref{lem:Y} or the lemmas
presented in \cite {BLY}, we obtain
$$
\norm{\nabla(\widetilde{ u}+v_3 )}_{L^{\infty} (\frac
1{r_1}\Omega\setminus (\frac 1 {r_1 }D_1 \cup \frac 1 {r_1 } D_2))}
\leq C \norm {\nabla H}_{L^{\infty}(\Omega)} \left(\frac {d}{1+d}
\right) \cdot \left\{
\begin{array}{ll}
\frac 1 {|\delta\log \delta|}~&\mbox{if}~ n= 3\notag\\
\frac 1 {\delta}~&\mbox{if}~ n \geq 4\notag
\end{array}
\right.$$ and $$\norm {\nabla v_3}_{L^{\infty}(\partial (\frac 1
{r_1 } D_1 \cup \frac 1 {r_1 } D_2)) } \leq C  \frac 1 d \norm
{\nabla H}_{L^{\infty}(\Omega)} \left( \frac {d}{1+d} { + 1}
\right).$$ Two bounds above yield
$$\norm {\nabla \widetilde{u}}_{L^{\infty}(\partial (\frac 1 {r_1 }D_1 \cup \frac 1 {r_1 } D_2)) } \leq
C \norm {\nabla H}_{L^{\infty}(\Omega)} \left(\frac {d}{1 + d}
\right)  \cdot \left\{
\begin{array}{ll}
\frac 1 {|\delta\log \delta|}~&\mbox{if}~ n= 3\notag\\
\frac 1 {\delta}~&\mbox{if}~ n \geq 4\notag
\end{array}
\right.$$ By the harmonicity of $ \widetilde{u} - \widetilde{H}$ in
$\mathbb{R}^n \setminus (\frac 1 {r_1} \Omega)$, we have
\begin{align*}\norm {\nabla( u - H )}_{L^{\infty}( \mathbb{R}^n
\setminus (D_1 \cup D_2)) } &= \norm {\nabla( \widetilde{u} -
\widetilde{H} )}_{L^{\infty}( \mathbb{R}^n \setminus (\frac 1 {r_1}
D_1 \cup \frac 1 {r_1} D_2)) }\\& \leq C \norm {\nabla
H}_{L^{\infty}(\Omega)} \left(\frac {r_1 r_2}{r_1 + r_2} \right)
\cdot \left\{
\begin{array}{ll}
\frac 1 {|\epsilon\log \epsilon|}~&\mbox{if}~ n= 3\notag\\
\frac 1 {\epsilon}~&\mbox{if}~ n \geq 4\notag
\end{array}
\right.\end{align*} Furthermore, this leads to
$$\norm {\nabla u }_{L^{\infty}( \Omega \setminus (D_1 \cup D_2)) }
\leq C \norm {\nabla H}_{L^{\infty}(\Omega)} \left(\frac {r_1
r_2}{r_1 + r_2} \right)  \cdot \left\{
\begin{array}{ll}
\frac 1 {|\epsilon\log \epsilon|}~&\mbox{if}~ n= 3\notag\\
\frac 1 {\epsilon}~&\mbox{if}~ n \geq 4\notag
\end{array}
\right.$$ Therefore, we have completed the proof. \qed

\section*{Acknowledgements}  The authors would like to express their gratitude to
\textrm{Professor YanYan Li} and his students Dr. Shiting Bao and
Biao Yin, who gave them a precious information on the optimal
blow-up rate in the higher dimensions and also encouraged the second
named author to check his result. The authors are also deeply
grateful to \textrm{Professor Hyeonbae Kang}, who suggested the
original subject studied in this paper and gave useful comments. The
second named author would also like to express his thanks to
\textrm{Professor David Kinderlehrer}, and gratefully acknowledges
his hospitality during the visiting period at Carnegie Mellon
University. The first author was supported by the Korea Research
Foundation Grant funded by the Korean Goverment (MOEHRD)
(KRF-2005-214-C00184).

\end{document}